# CONSTANT ANGLE RULED SURFACES IN EUCLIDEAN SPACES

Yusuf YAYLI [1], Evren ZIPLAR [2]

[1] *Department of Mathematics, Faculty of Science, University of Ankara, Tandoğan, Ankara, Turkey*

*yayli@science.ankara.edu.tr*

[2] *Department of Mathematics, Faculty of Science, University of Ankara, Tandoğan, Ankara, Turkey*

*evrenziplar@yahoo.com*

**Abstract.** In this paper, we study the special curves and ruled surfaces on helix hypersurface whose tangent planes make a constant angle with a fixed direction in Euclidean n-space $E^n$. Besides, we observe some special ruled surfaces in $\mathbb{IR}^n$ and give requirement of being developable of the ruled surface. Also, we investigate the helix surface generated by a plane curve in Euclidean 3-space $E^3$.

## 1. INTRODUCTION

Ruled surfaces are one of the most important topics of differential geometry. The surfaces were found by Gaspard Monge, who was a French mathematician and inventor of descriptive geometry. And, many geometers have investigated the many properties of these surfaces in [4,6,7].

Constant angle surfaces are considerable subject of geometry.There are so many types of these surfaces. Helix hypersurface is a kind of constant angle surfaces. An helix hypersurface in Euclidean n-space is a surface whose tangent planes make a constant angle with a fixed direction. The helix surfaces have been studied by Di Scala and Ruiz-Hernández in [2]. And, A.I.Nistor investigated certain constant angle surfaces constructed on curves in Euclidean 3-space $E^3$ [1].

One of the main purpose of this study is to observe the special curves and ruled surfaces on a helix hypersurface in Euclidean n-space $E^n$. Another purpose of this study is to obtain a helix surface by generated a plane curve in Euclidean 3-space $E^3$.

## 2. PRELIMINARIES

**Definition 2.1** Let $\alpha : I \subset \mathrm{IR} \to E^n$ be an arbitrary curve in $E^n$. Recall that the curve $\alpha$ is said to be of unit speed ( or parametrized by the arc-length function $s$) if $\langle \alpha'(s), \alpha'(s) \rangle = 1$, where $\langle , \rangle$ is the standart scalar product in the Euclidean space $E^n$ given by

$$\langle X, Y \rangle = \sum_{i=1}^{n} x_i y_i ,$$

for each $X = (x_1, x_2, ..., x_n), Y = (y_1, y_2, ..., y_n) \in E^n$.

Let $\{V_1(s), V_2(s), ..., V_n(s)\}$ be the moving frame along $\alpha$, where the vectors $V_i$ are mutually orthogonal vectors satisfying $\langle V_i, V_i \rangle = 1$. The Frenet equations for $\alpha$ are given by

$$\begin{bmatrix} V_1' \\ V_2' \\ V_3' \\ \vdots \\ V_{n-1}' \\ V_n' \end{bmatrix} = \begin{bmatrix} 0 & k_1 & 0 & 0 & \cdots & 0 & 0 \\ -k_1 & 0 & k_2 & 0 & \cdots & 0 & 0 \\ 0 & -k_2 & 0 & k_3 & \cdots & 0 & 0 \\ \vdots & \vdots & \vdots & \vdots & & \vdots & \vdots \\ 0 & 0 & 0 & 0 & \cdots & 0 & k_{n-1} \\ 0 & 0 & 0 & 0 & \cdots & -k_{n-1} & 0 \end{bmatrix} \begin{bmatrix} V_1 \\ V_2 \\ V_3 \\ \vdots \\ V_{n-1} \\ V_n \end{bmatrix}.$$

Recall that the functions $k_i(s)$ are called the $i$-th curvatures of $\alpha$ [3].

**Definition 2.2** Given a hypersurface $M \subset \mathrm{IR}^n$ and an unitary vector $d \neq 0$ in $\mathrm{IR}^n$, we say that $M$ is a helix with respect to the fixed direction $d$ if for each $q \in M$ the angle between $d$ and $T_q M$ is constant. Note that the above definition is equivalent to the fact that $\langle d, \xi \rangle$ is constant function along $M$, where $\xi$ is a normal vector field on $M$ [2].

**Theorem 2.1** Let $H \subset \mathrm{IR}^{n-1}$ be an orientable hypersurface in $\mathrm{IR}^{n-1}$ and let $N$ be an unitary normal vector field of $H$. Then,

$$f_\theta(x,s) = x + s(\sin(\theta)N(x) + \cos(\theta)d), \; f_\theta : H \times \mathrm{IR} \to \mathrm{IR}^n, \; \theta = \text{constant}$$

is a helix with respect to the fixed direction $d$ in $\mathrm{IR}^n$, where $x \in H$ and $s \in \mathrm{IR}$. Here $d$ is the vector $(0,0,...,1) \in \mathrm{IR}^n$ such that $d$ is orthogonal to $N$ and $H$ [2].

**Definition 2.3** Let $\alpha : \mathrm{I} \subset \mathrm{IR} \to E^n$ be a unit speed curve with nonzero curvatures $k_i$ $(i = 1,2,...,n)$ in $E^n$ and let $\{V_1, V_2, ... V_n\}$ denote the Frenet frame of the curve of $\alpha$. We call $\alpha$ a $V_n$-slant helix, if the $n$-th unit vector field $V_n$ makes a constant angle $\varphi$ with a fixed direction $\mathrm{X}$, that is,

$$\langle V_n, \mathrm{X} \rangle = \cos(\varphi), \; \varphi \neq \frac{\pi}{2}, \; \varphi = \text{constant}$$

along the curve, where $\mathrm{X}$ is unit vector field in $E^n$ [3].

## 3. THE SPECIAL CURVES ON THE HELIX HYPERSURFACES IN EUCLIDEAN n-SPACE $E^n$

**Theorem 3.1** Let $M$ be a helix hypersurface with the direction $d$ in $E^n$ and let $\alpha : \mathrm{I} \subset \mathrm{IR} \to M$ be a unit speed geodesic curve on $M$. Then, the curve $\alpha$ is a $V_2$-slant helix with the direction $d$ in $E^n$.

**Proof:** Let $\xi$ be a normal vector field on $M$. Since $M$ is a helix hypersurface with respect to $d$, $\langle d, \xi \rangle = $ constant. That is, the angle between $d$ and $\xi$ is constant on every point of the surface $M$. And, $\alpha''(s) = \lambda \xi |_{\alpha(s)}$ along the curve $\alpha$ since $\alpha$ is a geodesic curve on $M$. Moreover, by using the Frenet equation $\alpha''(s) = V_1' = k_1 V_2$, we obtain $\lambda\xi|_{\alpha(s)} = k_1 V_2$, where $k_1$ is first curvature of $\alpha$. Thus, from the last equation, by taking norms on both sides, we obtain $\xi = V_2$ or $\xi = -V_2$. So, $\langle d, V_2 \rangle$ is constant along the curve $\alpha$ since $\langle d, \xi \rangle = $ constant. In other words, the angle between $d$ and $V_2$ is constant along the curve $\alpha$. Consequently, the curve $\alpha$ is a $V_2$-slant helix with the direction $d$ in $E^n$.

**Corollary 3.1** For $n = 3$, the following Theorem obtained.
Theorem: Let $\alpha : I \subset \mathrm{IR} \to M$ be a curve on a constant angle surface $M$ with unit normal $N$ and the fixed direction $k$. If a curve $\alpha$ on $M$ is a geodesic, then $\alpha$ is a slant helix with the axis $k$ in $E^3$ (see [5]).

**Theorem 3.2** Let $M$ be a helix hypersurface in $E^n$ and let $\alpha : I \subset \mathrm{IR} \to M$ a unit speed curve on $M$. If the $n$-th unit vector field $V_n$ of $\alpha$ equals to $\xi$ or $-\xi$, where $\xi$ is a normal vector field on $M$, then $\alpha$ is a $V_n$-slant helix with the direction $d$ in $E^n$.

**Proof:** Let $d \neq 0 \in E^n$ be a fixed direciton of the helix hypersurface $M$. Since $M$ is a helix hypersurface with respect to $d$, $\langle d, \xi \rangle = $ constant. That is, the angle between $d$ and $\xi$ is constant on every point of the surface $M$. Let the $n$-th unit vector field $V_n$ of $\alpha$ be equals to $\xi$ or $-\xi$. Then $\langle d, V_n \rangle$ is constant along the curve $\alpha$ since $\langle d, \xi \rangle = $ constant. That is, the angle between $d$ and $V_n$ is constant along the curve $\alpha$. Finally, the curve $\alpha$ is a $V_n$-slant helix in $E^n$.

**Theorem 3.3** Let $M$ be a helix hypersurface with the direction $d$ in $E^n$ and let $\alpha : I \subset \mathrm{IR} \to M$ ($\alpha(t) \in M$, $t \in I$) be a curve on the surface $M$. If $\alpha$ is a line of curvature on $M$, then $d \in Sp\{T\}^\perp$ along the curve $\alpha$, where $T$ is tangent vector field of $\alpha$.

**Proof:** Since $M$ is a helix hypersurface with the direction $d$,

$$\langle N \circ \alpha, d \rangle = \text{constant}$$

along the curve $\alpha$, where $N$ is the normal vector field of $M$. If we are taking the derivative in each part of the equality with respect to $t$, we obtain:

$$\langle (N \circ \alpha)', d \rangle = 0.$$

Since $\alpha$ is a line of curvature on $M$, $(N \circ \alpha)' = S(T) = \lambda T$, where $S$ is the shape operator of the surface $M$. So, we have

$$\langle T, d \rangle = 0.$$

Finally, $d \in Sp\{T\}^{\perp}$ along the curve $\alpha$.

## 4. THE RULED SURFACES IN $\mathrm{IR}^n$

**Definition 4.1** Let $H \subset \mathrm{IR}^{n-1}$ be a orientable hypersurface in $\mathrm{IR}^{n-1}$ and let $\beta$ be a curve on $H$, where

$$\beta : I \subset \mathrm{IR} \to H \subset \mathrm{IR}^{n-1}.$$
$$t \to \beta(t)$$

Then,

$$\Phi(t, s) = \beta(t) + s(\sin(\theta) N(\beta(t)) + \cos(\theta) d)$$

is a ruled surface with dimension 2 on $f_\theta$ in $\mathrm{IR}^n$ ($f_\theta$ was defined in Theorem 2.1), where $\theta$ constant, $N$ is a unitary normal vector field of $H$ and $d$ is constant vector as defined in Theorem 2.1. The surface $\Phi$ will be called the ruled surface generated by the curve $\beta$.

**Theorem 4.1** The ruled surface $\Phi(t, s)$ defined above is developable if and only if the curve $\beta$ is a line of curvature on the surface $H$.

**Proof:** We assume that $\beta$ is a line of curvature on $H$.
Let consider the surface $\Phi(t, s) = \beta(t) + s(\sin(\theta) N(\beta(t)) + \cos(\theta) d)$ with rulings $X(t) = \sin(\theta) N(\beta(t)) + \cos(\theta) d$ and directrix $\beta$. If we are taking the partial derivative in each part of the equality with respect to $t$, we obtain:

$$\Phi_t = \frac{\partial \Phi}{\partial t} = \beta' + (s\sin(\theta))(N \circ \beta)' \tag{1}$$

And, $S(T) = \lambda T$ since $\beta$ is a line of curvature on $H$, where $T$ is tangent vector field of $\beta$ and $S$ is the shape operator of the surface $H$. Besides, $(N \circ \beta)' = \frac{dN}{dt} = S(T)$. Therefore, $(N \circ \beta)' = \lambda T$ and by using (1), we obtain the equality

$$\Phi_t = (1 + \lambda.s.\sin(\theta))\beta' = (1 + \lambda.s.\sin(\theta))T.$$

Hence, the system $\{\Phi_t, T\}$ is linear dependent. And, we know that a tangent plane along a ruling is spanned by $\Phi_s = X(t)$ and $\Phi_t$. Finally, the tangent planes are parallel along the ruling $\sin(\theta)N(\beta(t)) + \cos(\theta)d$ passing from the point $\beta(t)$. That is, the surface $\Phi(t,s)$ is developable.

We assume that the ruled surface $\Phi(t,s)$ is developable. Then, the system $\{\Phi_t, T\}$ is linear dependent. So, from the equality

$$\Phi_t = \frac{\partial \Phi}{\partial t} = \beta' + (s\sin(\theta))(N \circ \beta)' = T + (s\sin(\theta))(N \circ \beta)',$$

we get $(N \circ \beta)' = \lambda T$. Therefore $(N \circ \beta)' = S(T) = \lambda T$, where $S$ is the shape operator of the surface $H$. That is, $\beta$ is a line of curvature on $H$.

**Corollary 4.1** Let $H$ be the hypersphere

$$S^{n-2} = \left\{ x = (x_1, x_2, ..., x_{n-1}) : f(x) = \sum_{i=1}^{n-1} x_i^2 = 1, \left\| \vec{\nabla} f \right\| \neq 0 \right\} \subset \mathrm{IR}^{n-1}.$$

Let $\beta$ be a curve on $H = S^{n-2}$ where

$$\beta : I \subset \mathrm{IR} \to H = S^{n-2} \subset \mathrm{IR}^{n-1}$$
$$t \to \beta(t)$$

Then, the ruled surface $\Phi(t,s) = \beta(t) + s(\sin(\theta)\beta(t) + \cos(\theta) d) \subset \mathrm{IR}^n$ is always developable from Theorem 4.1. Because, each curve on the hypersphere $S^{n-2}$ is a line of curvature.

## 5. HELIX SURFACES GENERATED BY A PLANE CURVE IN EUCLIDEAN 3-SPACE $E^3$

Let
$$\alpha : I \subset \mathrm{IR} \to E^3$$
$$u \to \alpha(u)$$

be a plane curve in Euclidean 3-space $E^3$. And, we denote the tangent, principal normal, the binormal of $\alpha$ by $V_1 = T$, $V_2$ and $V_3 = B$, respectively. Note that binormal of a plane curve in $E^3$ is constant.

**Definition 5.1** We can obtain a ruled surface by using the plane curve $\alpha$ such that
$$\phi : U \subset E^2 \to E^3$$
$$(u,v) \to \phi(u,v) = \alpha(u) + v(\sin(\theta)V_2(u) + \cos(\theta)B).$$

The ruled surface will be called as the surface generated by the curve $\alpha$.

**Theorem 5.1** The ruled surface
$$\phi : U \subset E^2 \to E^3$$
$$(u,v) \to \phi(u,v) = \alpha(u) + v(\sin(\theta)V_2(u) + \cos(\theta)B)$$

is a helix surface with the direction $B$ in $E^3$, where $\theta$ is constant, $\alpha$ is a plane curve and $B$ is a constant vector which is perpendicular to the plane of the curve $\alpha$.

**Proof:** We want to show that $\langle Z, B \rangle$ is a constant function along $\phi$, where $Z$ is a normal vector field of $\phi$.

First, we are going to find a normal vector field $Z$. To do this, we will compute the partial derivatives of $\phi$ with respect to $u$ and $v$. Note that $k_2 = 0$ for the curve $\alpha$ since $\alpha$ is a planer curve in $E^3$.

$$\phi_u(u,v) = (1 - k_1.v.\sin(\theta))T \text{ and } \phi_v(u,v) = \sin(\theta)V_2 + \cos(\theta)B. \quad (2)$$

Using the equalities in (2), a normal to the surface $\phi$ is given by
$$Z = \frac{\phi_u \times \phi_v}{\|\phi_u \times \phi_v\|} = -\cos(\theta)V_2 + \sin(\theta)B.$$

So, we have $\langle Z, B \rangle = \sin(\theta) = const.$ Finally, $\phi$ is a helix surface with the direction $B$ in $E^3$.

This completes the proof.

**Corollary 5.1** The helix surface

$$\phi : U \subset E^2 \to E^3$$
$$(u,v) \to \phi(u,v) = \alpha(u) + v(\sin(\theta)V_2(u) + \cos(\theta)B).$$

is always developable.

**Proof:** We know that If $\det(T, X, X') = 0$, where $X = \sin(\theta)V_2 + \cos(\theta)B$ and $T$ tangent of $\alpha$, then $\phi$ is developable.

So, we will compute $\det(T, X, X')$:

$$T = \alpha'$$
$$X = \sin(\theta)V_2 + \cos(\theta)B$$
$$X' = -k_1 . \sin(\theta).T$$

and so, we have

$$\det(T, X, X') = 0.$$

This completes the proof.

**Theorem 5.2** Let

$$\alpha : I \subset \mathrm{IR} \to E^3$$
$$s \to \alpha(s)$$

be a plane curve (not a straight line) with unit speed in Euclidean 3-space $E^3$. We consider the helix surface (generated by the curve $\alpha(s)$)

$$\phi : U \subset E^2 \to E^3$$
$$(s, v) \to \phi(s, v) = \alpha(s) + v(\sin(\theta)V_2(s) + \cos(\theta)B).$$

Then, the Gauss curvature of $\phi$ is zero, and the mean curvature of $\phi$:

$$H = \frac{1}{2} \cdot \frac{k_1 . \cos(\theta)}{1 - k_1 . v . \sin(\theta)},$$

where $k_1$ is the first curvature of the curve $\alpha$.

**Proof:** From corollary 5.1, the surface $\phi$ is developable. So, the Gauss curvature of $\phi$ is zero.

Now, we are going to prove that $H = \frac{1}{2} \cdot \frac{k_1.\cos(\theta)}{1-k_1.v.\sin(\theta)}$.

The system $\{x_1, x_2\}$ is an orthonormal basis for the tangent space of $\phi$ at the point $\phi(s,v)$, where $x_1 = \frac{\phi_s}{\|\phi_s\|}$ and $x_2 = \phi_v$ ($\phi_s$ is partial derivative of $\phi$ with respect to $s$ and $\phi_v$ is partial derivative of $\phi$ with respect to $v$). Recall that a normal vector field of $\phi$ is $Z(s,v) = -\cos(\theta)V_2(s) + \sin(\theta)B$ by Theorem 5.1. And, we know that the mean curvature of $\phi$ at a point $\phi(s,v)$:

$$H(\phi(s,v)) = \frac{1}{2}\sum_{i=1}^{2}\langle S(x_i), x_i\rangle,$$

where $S$ is the shape operator of $\phi$.

So firstly, we will compute $S(x_1)$ and $S(x_2)$:

$$S(x_1) = D_{x_1}Z = D_{\frac{\phi_s}{\|\phi_s\|}}Z = \frac{1}{\|\phi_s\|}D_{\phi_s}Z = \frac{1}{\|\phi_s\|}\frac{dZ}{ds} = \left(\frac{k_1.\cos(\theta)}{|1-k_1.v.\sin(\theta)|}\right)T$$

and

$$S(x_2) = D_{x_2}Z = D_{\phi_v}Z = \frac{dZ}{dv} = 0,$$

where $D$ is standard covariant derivative in $E^3$.

Therefore, we have

$$\langle S(x_1), x_1\rangle = \frac{k_1.\cos(\theta)}{1-k_1.v.\sin(\theta)} \text{ and } \langle S(x_2), x_2\rangle = 0.$$

Finally,

$$H = \frac{1}{2}\sum_{i=1}^{2}\langle S(x_i), x_i\rangle = \frac{1}{2}\cdot\frac{k_1.\cos(\theta)}{1-k_1.v.\sin(\theta)},$$

where $1-k_1.v.\sin(\theta) \neq 0$.

This completes the proof.

**Corollary 5.2** The surface $\phi$ defined above is minimal if and only if $\theta = \pi/2$ where $1-k_1.v.\sin(\theta) \neq 0$. In that case (whenever $\theta = \pi/2$), the surface $\phi$ is a plane.

**Example 5.1** Let the curve $\alpha(u)$ be a plane curve parametrized by the vector function

$$\alpha(u) = \left(\frac{3}{5}\sin(u), 1+\cos(u), \frac{4}{5}\sin(u)\right), \ u \in [0, 5\pi].$$

Then,

$$V_2 = \left(-\frac{3}{5}\sin(u), -\cos(u), -\frac{4}{5}\sin(u)\right)$$

$$B = \left(\frac{4}{5}, 0, -\frac{3}{5}\right)$$

where $V_2$ is the principal normal and $B$ is the binormal of $\alpha$, respectively.

So, If we choose $\theta = \pi/6$ and $v \in [0, \pi]$, the helix surface generated by the curve $\alpha(u)$ has the parametric representation:

$$x = \left(\frac{3}{5} - \frac{3v}{10}\right)\sin(u) + \frac{2\sqrt{3}}{5}v$$

$$y = \left(1 - \frac{v}{2}\right)\cos(u) + 1$$

$$z = \left(\frac{4}{5} - \frac{2v}{5}\right)\sin(u) - \frac{3\sqrt{3}v}{10}.$$

And, the surface generated by the curve $\alpha$ is shown the following Figure.

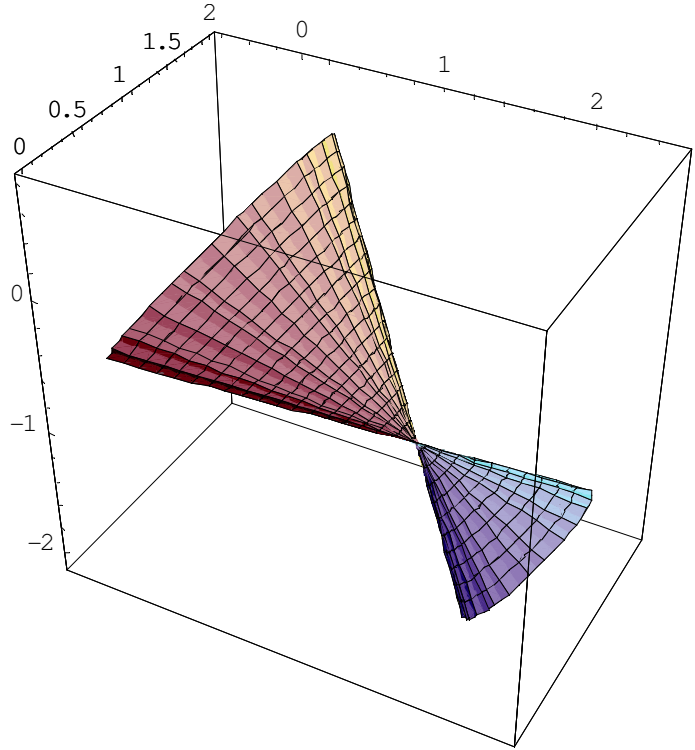

**REFERENCES**

[1] A.I. Nistor., 2009, Certain constant angle surfaces constructed on curves, arXiv:0904.1475v1 [math.DG]

[2] Di Scala, A.J., Ruiz-Hernández, G., 2009, Helix submanifolds of Euclidean spaces, Monatsh Math 157: 205-215.

[3] Gök, I., Camcı, Ç., Hacısalihoğlu, H.H., 2009, $V_n$-slant helices in Euclidean $n$-space $E^n$, Math. Commun., Vol. 14, No. 2, pp. 317-329.

[4] Holditch, A., Lady's and gentleman's diary for year 1858.

[5] Özkaldı, S., Yaylı, Y., 2011, Constant angle surfaces and curves in $E^3$, International electronic journal of geometry., Volume 4 No. 1 pp. 70-78.

[6] Sarıoğlugil, A., Tutar, A., 2007, On ruled surface in Euclidean space $E^3$, Int. J. Contemp. Math. Sci., Vol. 2, no. 1, pp. 1-11.

[7] Steiner, J., Ges. werke, Berlin, 1881-1882.